\documentclass[12pt]{article}
\usepackage{titling}
\usepackage[utf8]{inputenc}
\usepackage{amsmath}
\usepackage{amssymb} 
\usepackage{stmaryrd} 
\usepackage{mathtools} 
\usepackage{amsmath,amsthm,nicefrac}
\usepackage{hyperref} 

\usepackage{tikz}
\usepackage{tkz-euclide}
 \usepackage{cancel}
\usepackage[french]{babel}
\usepackage[T1]{fontenc}  
\usepackage{csquotes}
\usepackage{tcolorbox}
\usepackage{float} 
\usepackage[backend=biber,maxbibnames = 4,sorting=none]{biblatex} 
\usepackage{listings} 
\lstnewenvironment{pycode}{%
	\lstset{literate={é}{{\'e}}1{è}{{\`e}}1{ê}{{\^e}}1{ë}{{\"e}}1{à}{{\`a}}1{â}{{\^a}}1{ä}{{\"a}}1{á}{{\'a}}1{ï}{{\"i}}1{î}{{\^i}}1{ì}{{\`i}}1{í}{{\'i}}1{ù}{{\`u}}1{û}{{\^u}}1{ü}{{\"u}}1{ú}{{\'u}}1{ô}{{\^o}}1{ö}{{\"o}}1{ó}{{\'o}}1,numberstyle=\fontfamily{lmtt}\tiny,basicstyle=\fontfamily{lmtt}\small,language=python,numbers=left,%
		morekeywords={}%
}}{}

\usepackage{url}
\title{All of Monty}

\usetikzlibrary{calc,arrows.meta}

\usepackage{pgfplots}

\usepackage{mathpazo}
\usepackage{amsmath}

\newcommand{\pp}{\mathbb P}

\usepackage{pifont}

\author{Karim \textsc{Zayana}\up{1,2} et Ivan \textsc{Boyer}\up{1}}

\date{\small\up{1} Ministère de l'\'Education nationale, Paris\\
\up{2} LTCI, Télécom Paris, Institut Polytechnique de Paris}

\definecolor{myblue}{RGB}{79,129,189}
\definecolor{mygreen}{RGB}{0,176,80}
\begin{document}
\maketitle

Tout a déjà été dit sur le problème de Monty Hall, au point d'attirer l'attention du programme de mathématiques intégré à l'enseignement scientifique de première \cite{progPremEnsSciMath}. L'objet de ce texte n'est donc pas d'en rajouter, mais plutôt d'y faire le tri tout en clarifiant quelques scénarios pédagogiques parus dans la littérature -- foisonnante -- déjà consacrée au sujet \cite{rosenhouse2009, rosenhouse2009bis}.

\section{L'origine du problème}
Le problème est inspiré d'une célèbre émission, <<~Let's Make a Deal~>>, diffusée successivement sur les chaînes américaines NBC puis ABC de 1963 à 1977. Son présentateur vedette, le canadien Monty Hall (de son vrai nom Maurice Halprin) y campe un aussi joyeux qu'habile bonimenteur. Circulant parmi le public, le verbe haut dans ses habits de lumière, il alpague le chaland. Et ouvre les paris. Facétieux, l'animateur remet des enveloppes cadeau différentes à deux invités... avant de leur souffler d'échanger... à leurs risques et périls ! Roublard, il offre $100$~\$ à une spectatrice : qu'elle les conserve... ou bien les dépense pour ce que cache une boîte qu'il lui tend aussitôt. Peut-être contient-elle une rivière de diamants... à moins qu'il ne s'y trouve qu'un petit grelot, <<~Who knows?~>>.  Madame hésite. Il renchérit : <<~Well I give you some \$$100$ more, then you keep the money and we end the game, all'right ?~>>. Cruel, il dévoile avec componction l'arrière du rideau d'entre les trois qu'aurait dû choisir un autre candidat. Le malheureux ne repartira qu'avec ses rêves de millionnaire et pour seul prix un sachet de bonbons. Et très théâtralement, la caméra de saisir les mines tantôt déconfites, tantôt pâmoisées de cette grande comédie.   

\bigskip
Ce show télévisé essaimera par la suite sous forme de reprises, plus ou moins heureuses, aux USA comme ailleurs. En France par exemple, son adaptation, le <<~Bigdil~>> (contraction de <<~Big~>> et <<~Deal~>>), battra des records d'audience dans les années 2000 sur la chaîne TF1. 

\bigskip
Mais c'est dès 1975 que les mathématiques s'emparèrent du concept. Tout commence par un courrier des lecteurs publié dans une revue de statistiques \cite{Selvin1975}. L'auteur, de l'université de Berkeley, y signe un petit sketch entre les protagonistes du jeu. Parodique, le script pousse Monty Hall (rebaptisé <<~Monte Hall~>> pour la circonstance) à une escalade burlesque qui s'achève sur un dilemme absolument diabolique. Bien que l'article se fît à l'époque discret, ne sortant pas d'une communauté de scientifiques avertis, nous le reproduisons ici à l'identique et dans son intégralité tant il fut fondateur.    

\medskip
{\tt
\begin{center}
A PROBLEM IN PROBABILITY
\end{center}
It is <<~Let's Make a Deal~>> - a famous TV show starring Monte\footnote{La transformation du nom est d'origine (et sans doute volontaire).} Hall.

\medskip
\noindent\begin{minipage}[t]{0.3\linewidth}
\vspace{0cm}Monte Hall:
\end{minipage}
\hfill
\begin{minipage}[t]{0.7\linewidth}
\vspace{0cm}One of the three boxes labeled A, B and C contains the keys to that new 1975 Lincoln Continental\footnote{Modèle sorti des usines Ford. Il fut commercialisé en versions berline, cabriolet, limousine, présidentielle (dont la décapotable à bord de laquelle John F. Kennedy roulait quand il fut assassiné à Dallas en 1963).}. The other two are empty. If you choose the box containing the keys, you win the car.\medskip\vfill
\end{minipage}

\noindent\begin{minipage}[t]{0.3\linewidth}
\vspace{0cm}Contestant: 
\end{minipage}
\hfill
\begin{minipage}[t]{0.7\linewidth}
\vspace{0cm}Gasp!\medskip\vfill
\end{minipage}

\noindent\begin{minipage}[t]{0.3\linewidth}
\vspace{0cm}Monte Hall:
\end{minipage}
\hfill
\begin{minipage}[t]{0.7\linewidth}
\vspace{0cm}Select one of these boxes.\medskip\vfill
\end{minipage}

\noindent\begin{minipage}[t]{0.3\linewidth}
\vspace{0cm}Contestant: 
\end{minipage}
\hfill
\begin{minipage}[t]{0.7\linewidth}
\vspace{0cm}I'll take box B.\medskip\vfill
\end{minipage}

\noindent\begin{minipage}[t]{0.3\linewidth}
\vspace{0cm}Monte Hall:
\end{minipage}
\hfill
\begin{minipage}[t]{0.7\linewidth}
\vspace{0cm}Now box A and box C are on the table and here is box B (contestant grips box B tightly). It is possible the car keys are in that box! I'll give you \$100 for the box.\medskip\vfill
\end{minipage}

\noindent\begin{minipage}[t]{0.3\linewidth}
\vspace{0cm}Contestant: 
\end{minipage}
\hfill
\begin{minipage}[t]{0.7\linewidth}
\vspace{0cm}No, thank you.\medskip\vfill
\end{minipage}

\noindent\begin{minipage}[t]{0.3\linewidth}
\vspace{0cm}Monte Hall:
\end{minipage}
\hfill
\begin{minipage}[t]{0.7\linewidth}
\vspace{0cm}How about \$200?\medskip\vfill
\end{minipage}

\noindent\begin{minipage}[t]{0.3\linewidth}
\vspace{0cm}Contestant: 
\end{minipage}
\hfill
\begin{minipage}[t]{0.7\linewidth}
\vspace{0cm}No!\medskip\vfill
\end{minipage}

\noindent\begin{minipage}[t]{0.3\linewidth}
\vspace{0cm}Audience: 
\end{minipage}
\hfill
\begin{minipage}[t]{0.7\linewidth}
\vspace{0cm}No!\medskip\vfill
\end{minipage}

\noindent\begin{minipage}[t]{0.3\linewidth}
\vspace{0cm}Monte Hall:
\end{minipage}
\hfill
\begin{minipage}[t]{0.7\linewidth}
\vspace{0cm}Remember that the probability of your box containing the keys to the car is $\nicefrac{1}{3}$ and the probability of your box being empty is $\nicefrac{2}{3}$. I'll give you \$500.\medskip\vfill
\end{minipage}

\noindent\begin{minipage}[t]{0.3\linewidth}
\vspace{0cm}Audience: 
\end{minipage}
\hfill
\begin{minipage}[t]{0.7\linewidth}
\vspace{0cm}No!!\medskip\vfill
\end{minipage}

\noindent\begin{minipage}[t]{0.3\linewidth}
\vspace{0cm}Contestant: 
\end{minipage}
\hfill
\begin{minipage}[t]{0.7\linewidth}
\vspace{0cm}No, I think I'll keep this box.\medskip\vfill
\end{minipage}

\noindent\begin{minipage}[t]{0.3\linewidth}
\vspace{0cm}Monte Hall:
\end{minipage}
\hfill
\begin{minipage}[t]{0.7\linewidth}
\vspace{0cm}I'll do you a favor and open one of the remaining boxes on the table (he opens box A). It's empty! (Audience: applause). Now either box C or your box B contains the car keys. Since there are two boxes left, the probability of your box containing the keys is
now $\nicefrac{1}{2}$. I'll give you \$1000 cash for your box.\medskip\vfill
\end{minipage}

\begin{center}
WAIT!!!!
\end{center}

Is Monte right? The contestant knows that at least one of the boxes on the table is empty. He now knows that it was box A. Does this
knowledge change his probability of having the box containing the keys from $\nicefrac{1}{3}$ to $\nicefrac{1}{2}$? One of the boxes on the table has to be empty. Has Monte done the contestant a favor by showing him which of the two boxes was empty? Is the probability of winning the car $\nicefrac{1}{2}$ or $\nicefrac{1}{3}$?

\noindent\begin{minipage}[t]{0.3\linewidth}
\vspace{0cm}Contestant:
\end{minipage}
\hfill
\begin{minipage}[t]{0.7\linewidth}
\vspace{0cm}I'll trade you my box B for the box C on the table.\medskip\vfill
\end{minipage}

\noindent\begin{minipage}[t]{0.3\linewidth}
\vspace{0cm}Monte Hall:
\end{minipage}
\hfill
\begin{minipage}[t]{0.7\linewidth}
\vspace{0cm}That's weird!!\medskip\vfill
\end{minipage}
\medskip
\noindent HINT:The contestant knows what he is doing!
}

Le sujet passa donc inaperçu jusqu'en en 1990, pour rejaillir inopinément dans les colonnes d'un magazine en vogue, Parade, à la rubrique <<~Ask Marilyn~>> (un supplément du dimanche, qui existe toujours) sous la plume de sa chroniqueuse, Marilyn Vos Savant \cite{Marilyn1990}. L'enjeu demeure le même, sauf que les boîtes y sont remplacées par des portes, les clés de voiture par la voiture et les vides par des chèvres. Très vite, le lectorat s'enflamme, comme le relatent le New York Times et le Statistician American d'alors \cite{Tierney1991, Morgan1991}. Les lettres passionnées affluent par milliers à la rédaction. On s'emporte contre l'automobile. On défend la cause animale. On se déchaîne sur la solution, pourtant correcte, fournie par le journal. Appelés à la rescousse, les experts en tous poils s'affrontent à grandes querelles de nombres.  Monty Hall lui-même reprendra du service, distillant \href{https://www.youtube.com/watch?v=c1BSkquWkDo}{interviews}, mémoires et confidences. C'est le buzz !

\bigskip
Le pastiche de 1975, ou sa déclinaison de 1990, ne seront pas joués dans l'émission d'origine. Le problème de Monty Hall n'est donc jamais qu'une expérience de pensée. Son dénouement, nous le verrons, dépend pour beaucoup de l'humeur du présentateur, et de la connaissance supposée que nous en avons, au moment d'ouvrir la boîte A. Une fois ces informations traduites en hypothèses mathématiques, nous constaterons que la probabilité recherchée balance, en effet, entre plusieurs valeurs : non pas deux seulement comme la fin du dialogue le laisse d'ailleurs imaginer, mais une infinité. 

\section{Le traitement du problème}
Il existe une multitude de façons d'aborder le problème, lesquelles mènent à d'innombrables variantes -- occasionnellement tirées par les cheveux. Nous nous limiterons à trois, mais nous les ferons toutes commencer de la même manière : le concurrent (son titre dans la distribution des rôles) montre la boîte ici baptisée B ; il la prend ; on pose à sa gauche la boîte ici baptisée A, à sa droite la boîte ici baptisée C. Quitte à ré-étiqueter les boîtes d'une manche à l'autre, nous tiendrons cet ordre pour immuable afin qu'il coïncide toujours avec la réalité physique de la scène. 

\bigskip

\bigskip
\textbf{2.1\quad Je (moi, Monty) sais que tu (toi, candidat) ne sais pas et tu  sais que je  sais}. Dans cette première version, qui est à la fois la plus courante et la plus naturelle, nous considérons que :
\begin{itemize}
\item le concurrent ne sait pas où se trouvent les clés de contact de la Lincoln Continental. Il choisit donc la boîte B par hasard;
\item Monty Hall sait en revanche où elles se trouvent. D'entre les deux boîtes restantes, A et C, il connaît celle(s) qui est(sont) vide(s). Il ouvre délibérément celle qui l'est, et indifféremment l'une ou l'autre quand il a le choix (c'est-à-dire quand les clés sont dans la boîte B désignée par le concurrent). 
\end{itemize}

\bigskip
Si le candidat trouve d'emblée la bonne boîte, en changer le fait perdre. Sinon, en changer le fait gagner avec certitude. Le candidat n'avait qu'une chance sur trois de trouver la bonne boîte. Il avait, en contrepoint, deux chances sur trois de se tromper. La stratégie consistant à changer de boîte inverse sa probabilité de gagner de $\nicefrac{1}{3}$ à $\nicefrac{2}{3}$. L'indice (<<~HINT~>>) qui conclut l'article de 1975 était bel et bien pertinent : l'intérêt du joueur est d'opter pour la boîte C lorsque Monty Hall ouvre la A. 

\begin{figure}[H]
\centering
\includegraphics[width =0.8\linewidth]{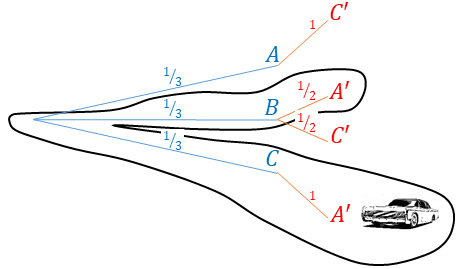}
\caption{Arbre de probabilités (les branches dites <<~mortes~>>, terminant par $B'$, ou encore $C'$ si issue de $C$, ne sont pas représentées car de valence nulle). Le premier choix du candidat se porte invariablement sur la boîte B. Monty Hall ouvre ici la boîte A. Sous cette condition, la probabilité que la boîte C soit gagnante calcule le <<~poids~>> du chemin obligeant Monty à ouvrir la boîte A rapporté au <<~poids~>> total des chemins conduisant Monty à ouvrir la boîte A. C'est la formule de Bayes !}
\label{Modele1}
\end{figure}

\bigskip
Nous pourrions nous satisfaire de ce raisonnement : élémentaire, il n'en est pas moins juste. Mais tant dans un but pédagogique qu'en prévision des deux modèles suivants, formalisons-le grâce au vocabulaire ad-hoc tout en l'étayant d'un arbre des possibles, figure \ref{Modele1}. À cet effet, notons 
\medskip
\begin{itemize}
\item l'événement $A$ : <<~La boîte A contient les clés~>>;
\item l'événement $B$ : <<~La boîte B contient les clés~>>;
\item l'événement $C$ : <<~La boîte C contient les clés~>>;
\item l'événement $A'$ : <<~Monty Hall ouvre la boîte A~>>;
\item l'événement $B'$ : <<~Monty Hall ouvre la boîte B~>> (à exclure);
\item l'événement $C'$ : <<~Monty Hall ouvre la boîte C~>>.
\end{itemize}  

\bigskip
La probabilité qu'a le candidat de remporter la voiture en changeant de boîte sachant que Monty Hall ouvre la A est $\mathbb{P}(C \mid A')$. Elle vaut aussi, d'après la formule de Bayes\footnote{Thomas \textsc{Bayes}, pasteur britannique du XVIII\ieme{} siècle. Mathématicien amateur, <<~sa~>> formule fut publiée à titre posthume.},
\begin{equation}
\mathbb{P}(C \mid A') = \frac{\mathbb{P}(A'\mid C)}{\mathbb{P}(A')}\mathbb{P}(C)
\label{Bayes}
\end{equation}
Dans sa philosophie, cette formule corrige \emph{a posteriori} la probabilité \emph{a priori} de l'événement $C$ à l'aide d'une information remontée à l'observateur (ici le joueur) : l'événement $A'$ \cite{zayana2020Cote}. La formule des probabilités totales transforme (\ref{Bayes}) en 
\begin{align*}
\mathbb{P}(C \mid A') &= \frac{\mathbb{P}(A'\mid C)\,\mathbb{P}(C)}{\mathbb{P}(A'\mid C)\,\mathbb{P}(C)+\mathbb{P}(A'\mid B)\,\mathbb{P}(B)}\\
&= \frac{1\times \nicefrac{1}{3}}{1\times \nicefrac{1}{3} + \nicefrac{1}{2}\times\nicefrac{1}{3}}\\
&= \nicefrac{2}{3}.
\end{align*}

Nous pouvions également nous en remettre à la complétude des événements $A$, $B$, $C$ relativement à la probabilité (conditionnelle) $\mathbb{P}(\bullet \mid A')$. En effet, puisque $\mathbb{P}(A \mid A')=0$ et que $\mathbb{P}(B \mid A')$ demeure inchangée à la valeur $\mathbb{P}(B)=\nicefrac{1}{3}$, l'égalité 
\begin{equation*}
\mathbb{P}(A \mid A')+\mathbb{P}(B \mid A')+\mathbb{P}(C \mid A')=1.
\end{equation*}
retourne directement
\begin{equation*}
\mathbb{P}(C \mid A')=\nicefrac{2}{3}.
\end{equation*}
Tout se passe comme si une partie des possibles, rendue impossible par la réalisation de l'événement $A'$, s'était transférée sur l'événement $C$.

\bigskip
Le candidat ne jouera le plus souvent qu'une unique fois. Mais ce qui intéresse la production, c'est l'espoir (techniquement, son espérance) de gain des candidats (qui défileront sur le plateau) quand ils optent pour une stratégie du changement systématique. Introduisons à cet effet la variable aléatoire $G$ égale à $1$ en cas de gain, à $0$ sinon. S'agissant d'une variable de Bernoulli\footnote{Jakob \textsc{Bernoulli} vécut en Suisse au XVII\ieme{} siècle et fut le premier d'une lignée de mathématiciens.}, $E(G) = \pp(G=1)$. Calculons, avec la formule des probabilités totales
\begin{align*}
\pp(G=1)&=\pp(G=1\mid A)\pp(A)+\pp(G=1\mid B)\pp(B)+\pp(G=1\mid C)\pp(C)\\
&= 1\times\nicefrac13+0\times\nicefrac13+1\times\nicefrac13=\nicefrac23.
\end{align*}
En effet, si les clés sont en A, Monty ouvre C, le candidat l'emporte en changeant  pour A ; si les clés sont en B, le candidat perd en changeant que Monty ouvre A ou C ; si les clés sont en C, Monty ouvre A, le candidat l'emporte en changeant pour C.

Dans l'autre scénario, où le candidat s'entête, cette espérance ne vaudrait que $\nicefrac13$. La première stratégie est donc payante.

\bigskip
S'il subsistait le moindre doute, voici une simulation informatique où sont répétées un million de fois les deux stratégies afin d'estimer, sur chacune, l'espérance de gain du candidat. Dans le second cas, celui où l'on garde toujours la boîte B, on comprend  sans même exécuter le programme que le résultat vaudra $\nicefrac13$, l'appel à la fonction \verb@choixMonty@ n'ayant pas incidence. La première stratégie conduit, elle, à une espérance  de  $\nicefrac23$ comme on l'a expliqué. Dans le script ci-après, la boîte A est notée $0$, la B, $1$ et la C, $2$. Le candidat s'arrête donc d'abord sur l'urne numérotée $1$, et Monty Hall choisit l'une de celles numérotées $0$ ou $2$. 

\begin{pycode}
import random as rd

def choixMonty(boiteVoiture):
    if boiteVoiture == 1:
        return rd.randint(0,1)*2  # 0 ou 2 au hasard
    return 2-boiteVoiture         # l'autre parmi 0 ou 2

def strategie1(): # on change de boîte
    boiteVoiture=rd.randint(0,2)
    boiteMonty=choixMonty(boiteVoiture)
    return 2-boiteMonty == boiteVoiture  
    # la dernière boîte est 0+1+2-boiteCandidat-boiteMonty
    # soit 2-boiteMonty

def strategie2(): # on garde la boîte
    boiteVoiture=rd.randint(0,2)
    boiteMonty=choixMonty(boiteVoiture) 
    # ici on voit que cette ligne ne sert à rien : 
    # le bonimenteur ne change pas les probabilités !
    return 1==boiteVoiture

def simulation1(N):
    gains=0
    for i in range(N):
        if strategie1():
            gains+=1
    return gains/N

print(simulation1(10**6),"proche de 2/3")

def simulation2(N):
    gains=0
    for i in range(N):
        if strategie2():
            gains+=1
    return gains/N

print(simulation2(10**6),"proche de 1/3")
\end{pycode}

\bigskip
\textbf{2.2\quad Je sais que tu ne sais pas et tu  sais que je  ne sais pas}. Dans cette deuxième version, nous considérons toujours que
\begin{itemize}
\item le concurrent ne sait pas où se trouvent les clés de contact de la Lincoln Continental. Il choisit donc la boîte B par hasard;
\item Monty Hall sait en revanche où elles se trouvent. Mais il est gaucher. Faisant face au candidat, il ouvre alors plus volontiers la boîte de droite (soit, la C) quand il a le choix. Cette tendance est mesurée par une probabilité $p$ proche de $1$ connue du joueur. En effet, en compétiteur coriace, ce dernier aura visionné toutes les archives. Cependant que Monty Hall, peu enclin aux chiffres et donc inconscient du biais qu'introduit sa routine, ne cherchât à la rectifier.
\end{itemize}

\begin{figure}[H]
\centering
\includegraphics[width =0.8\linewidth]{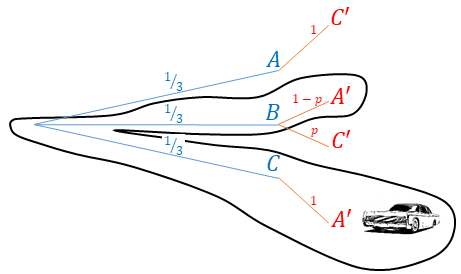}
\caption{Arbre de probabilités. Gaucher et faisant face à l'invité, Monty Hall marque une propension $p$ (inconsciente mais connue du joueur) à ouvrir la boîte C quand il a le choix.}
\label{Modele2}
\end{figure}

\medskip
L'arbre des possibles décrivant la situation est représenté en figure \ref{Modele2}. Il s'en suit
\begin{align*}
\mathbb{P}(C \mid A') &= \frac{1\times \nicefrac{1}{3}}{1\times \nicefrac{1}{3} + (1-p)\times\nicefrac{1}{3}}\\
&= \frac{1}{2-p}.
\end{align*}

Le modèle 1 correspond au modèle 2 dans le cas où Monty est ambidextre ($p=\nicefrac{1}{2}$). Et dans le cas extrême où $p=1$, $\mathbb{P}(C \mid A')$ grimpe à la valeur $1$. Logique : l'animateur ouvre la boîte A si, et seulement si, il y est contraint. De fait le gros lot est en C. 

Mais peut-on tirer avantage de cette manie en termes, cette fois, de gain moyen ? Les candidats seraient tentés d'adopter cette stratégie systématique : si Monty ouvre la boîte A (ce qui correspond à l'événement $A'$), changer d'office pour la boîte C semble une bonne option ; si, par contre, Monty ouvre la boîte C (ce qui correspond à l'événement $C'$), conserver la boîte $B$ paraît avisé. Un tel algorithme conduit à l'espérance suivante (en reprenant les notations de la section précédente)~:

\begin{align*}
\pp(G=1)&=\pp(G=1\mid A)\pp(A)+\pp(G=1\mid B)\pp(B)+\pp(G=1\mid C)\pp(C)\\
&= 0 \times \nicefrac13 + p \times \nicefrac13 + 1 \times \nicefrac13 \\
&= \tfrac{p+1}3.
\end{align*}
En effet, si les clés sont en A, Monty ouvre C, le candidat reste sur B et perd ; si les clés sont en B, Monty ouvre A (resp. C) avec la propension $1-p$ (resp. $p$) et le candidat change pour C (resp. reste en B) : il perd (resp. gagne) ; si les clés sont en C, Monty ouvre A, le candidat l'emporte en changeant pour C.

Pour $p$ voisin de $1$, l'espérance avoisine ainsi $\nicefrac23$. Nonobstant, on ne fait pas mieux que les $\nicefrac23$ calculés au paragraphe précédent : l'espérance n'a pas profité du penchant de Monty, elle en a même souffert ! De quoi rassurer la production...

\bigskip
\textbf{2.3\quad Je sais que tu ne sais pas et tu  ne sais pas que je  sais}. Dans cette troisième et dernière version, nous considérons toujours que
\begin{itemize}
\item le concurrent ne sait pas où se trouvent les clés de contact de la Lincoln Continental. Il choisit donc la boîte B par hasard;
\end{itemize}
Toutefois, 
\begin{itemize}
\item le concurrent ne sait pas si Monty choisit en connaissance de cause entre les boîtes A et C : la survenance de l'événement $A'$ n'induit plus implicitement celle du complémentaire $\overline{A}$. Il doit ainsi admettre que c'est fortuitement qu'ouvrant la boîte A celle-ci s'avérait vide. 
\end{itemize}

\begin{figure}[H]
\centering
\includegraphics[width =0.8\linewidth]{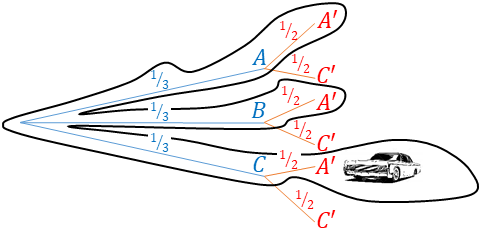}
\caption{Arbre de probabilités. Monty ouvre par hasard l'une des deux boîte A ou C. Le risque est payant (pour lui).}
\label{Modele3}
\end{figure}

\bigskip
Le candidat ne peut dès lors que se fonder sur l'arbre des possibles, au goût brouillé de la figure \ref{Modele3}. Détaillons :
\begin{align*}
\mathbb{P}(C \mid A' \cap \overline{A}) &=  \frac{\mathbb{P}(A' \cap \overline{A}\mid C) \, \mathbb{P}(C)}{\mathbb{P}(A' \cap \overline{A}\mid A) \, \mathbb{P}(A) + \mathbb{P}(A' \cap \overline{A}\mid B) \, \mathbb{P}(B) + \mathbb{P}(A' \cap \overline{A}\mid C) \, \mathbb{P}(C)}\\
&=\frac{\mathbb{P}(A' \mid C) \, \mathbb{P}(C)}{0 \times \mathbb{P}(A) + \mathbb{P}(A' \mid B) \, \mathbb{P}(B) + \mathbb{P}(A' \mid C) \, \mathbb{P}(C)}\\
&= \frac{\nicefrac{1}{2}\times \nicefrac{1}{3}}{0+\nicefrac{1}{2}\times \nicefrac{1}{3}+\nicefrac{1}{2}\times \nicefrac{1}{3}}\\
&=\nicefrac{1}{2}.
\end{align*}

\noindent C'est du fifty-fifty, Monty a réussi son coup !

\bigskip
Karim \textsc{Zayana} est inspecteur général de l'éducation, du sport et de la recherche (groupe des mathématiques) et professeur invité à l'institut polytechnique de Paris (Palaiseau) et Ivan \textsc{Boyer} est professeur de mathématiques en MPSI au lycée Champollion (Grenoble).
\printbibliography 
\end{document}